\documentclass{proc-l}

\issueinfo{00}
  {0}
  {Xxxxx}
  {2007}
\PII{S \ISSN}
\copyrightinfo{2007}{American Mathematical Society}
\pagespan{0}{0}
\commby{Joseph A. Ball}
\date{December 20, 2006}
\usepackage{hyperref}

\newtheorem{theorem}{Theorem}
\newtheorem{lemma}[theorem]{Lemma}

\newtheorem{remark}[theorem]{Remark}

\newcommand{\R}{{\mathbb R}}

\newcommand{\C}{{\mathbb C}}

\newcommand{\be}{\begin{equation}}
\newcommand{\ee}{\end{equation}}
\newcommand{\bea}{\begin{eqnarray}}
\newcommand{\eea}{\end{eqnarray}}

\newcommand{\ti}{\tilde}

\newcommand{\id}{{\rm 1\hspace{-0.6ex}l}}

\newcommand{\I}{\mathrm{i}}
\newcommand{\tr}{\mathrm{tr}}

\newcommand{\db}{\mathfrak{D}}
\DeclareMathOperator{\Ran}{Ran}

\newcommand{\sig}{\sigma}
\newcommand{\lam}{\lambda}
\newcommand{\gam}{\gamma}

\begin{document}

\title[Approximation of Isolated Eigenvalues]{On the Approximation of Isolated Eigenvalues
of Ordinary Differential Operators}

\author[G. Teschl]{Gerald Teschl}
\address{Faculty of Mathematics\\
Nordbergstrasse 15\\ 1090 Wien\\ Austria\\ and International Erwin Schr\"odinger
Institute for Mathematical Physics, Boltzmanngasse 9\\ 1090 Wien\\ Austria}
\email{\href{mailto:Gerald.Teschl@univie.ac.at}{Gerald.Teschl@univie.ac.at}}
\urladdr{\href{http://www.mat.univie.ac.at/~gerald/}{http://www.mat.univie.ac.at/\~{}gerald/}}

\thanks{{\it Research supported by the Austrian Science Fund (FWF) under Grant No.\ Y330}}

\keywords{Sturm--Liouville operators, Dirac operators, eigenvalues}
\subjclass[2000]{Primary 34L40, 34L16; Secondary  47N50, 34B20}

\begin{abstract}
We extend a result of Stolz and Weidmann on the approximation of isolated eigenvalues of singular
Sturm--Liouville and Dirac operators by the eigenvalues of regular operators.
\end{abstract}

\maketitle

\section{Introduction}

The approximation of isolated eigenvalues of singular ordinary differential operators
by the eigenvalues of regular operators is an important and well studied topic since
the latter ones can be computed numerically with arbitrary precision. See the recent
monograph by Zettl \cite{ze} or in particular the recent survey \cite{wd2} by Weidmann.

While the case of eigenvalues below the essential spectrum is well understood, the case
of eigenvalues in essential spectral gaps was only recently solved by Stolz and Weidmann
in \cite{sw} (see also \cite{sw2}).

Let $\tau$ be either a Sturm--Liouville expression
\be
\tau = \frac{1}{r} \Big(- \frac{d}{dx} p \frac{d}{dx} + q \Big),
\ee
or a Dirac system
\be
\tau = \frac{1}{r} \Big( \I \sig_2 \frac{d}{dx} + q \Big),
\ee
where $\sig_2=\begin{pmatrix} 0 & -\I\\ \I & 0\end{pmatrix}$ is the Pauli matrix, on the interval $I=(a,b)$.

As usual, we will assume  the coefficients $p^{-1},q,r$ are real-valued locally integrable
functions with $p,r>0$ in the Sturm-Liouville case and $q,r$ are real, symmetric $2\times2$
matrices with $r>0$ in the Dirac case.

An endpoint $a$ or $b$ is called regular if it is finite and the coefficients are integrable near
this endpoint. If both endpoints are regular, we will call $\tau$ regular.

Let $\db(\tau)$ be the maximal domain of definition of $\tau$. Then $\tau$ is self-adjoint
on $\db(\tau)$ if it is limit point (l.p.) near both $a$ and $b$. Otherwise we will impose an
additional boundary condition at every endpoint where $\tau$ is limit circle (l.c.). In this
way we obtain a self-adjoint operator $H$ associated with $\tau$.

A solution $\psi_a(z,x)$ ($\psi_b(z,x)$) of $\tau \psi = z\psi$ which is square integrable near
$a$ ($b$) and satisfies the boundary condition at $a$ ($b$) (if any) is called a Weyl solution.
Such a solution is unique up to a constant and it exists at least for $z\in\C\backslash\sig_{ess}(H)$.

Our aim is to approximate $H$ by regular operators $H_n$ obtained by restricting
$\tau$ to a finite interval $(a_n,b_n) \subseteq (a,b)$ (\cite[Chap.~14]{wdln}).
The case $a_n=a$ or $b_n=b$ is allowed.

Fix functions $u,v\in\db(\tau)$. Pick $a_n\downarrow a$,
$b_n\uparrow b$. Define $H_n$
\be \label{deftiHm}
H_n: \begin{array}[t]{lcl} \db(H_n) &\to& L^2((a_n,b_n); r\, dx) \\ f
&\mapsto& \tau f \end{array},
\ee
where
\be
\db(H_n) = \{ f \in \db(\tau_n)) | W_{a_n}(u,f) = W_{b_n}(v,f) =0 \}.
\ee

Then Stolz and Weidmann prove the following

\begin{theorem}[\cite{sw}]\label{thm:sw}
Define $H_n$ as above with $u=\psi_a(\lam_a)$ and $v=\psi_b(\lam_b)$ with
$\lam_a,\lam_b \in [\lam_0,\lam_1]$ (in particular, assume that the corresponding Weyl
solutions exist). Let $P_\Omega(H)$ be the spectral projection of $H$ corresponding
to the Borel set $\Omega \subseteq\R$.

If $\dim\Ran P_{(\lam_0,\lam_1)}(H)<\infty$, then the eigenvalues of $H$
in $(\lam_0,\lam_1)$ are exactly the limits of eigenvalues of $H_n$ which lie in
$(\lam_0,\lam_1)$. The corresponding (one-dimensional) eigenprojections converge
in norm.

If $\dim\Ran P_{(\lam_0,\lam_1)}(H)=\infty$, then the eigenvalues of $H_n$
accumulate in $(\lam_0,\lam_1)$ as $n\to\infty$.
\end{theorem}

\noindent
In fact, in \cite{sw} this result is only proven in the case $[\lam_0,\lam_1]\cap\sig_{ess}(H) = \emptyset$
and $\lam_0$, $\lam_1$ are not eigenvalues of $H$. Hence we will provide a proof of this
slightly generalized version below.

As pointed out in \cite{sw}, this result has of course one practical drawback: The
Weyl solutions used to generate the boundary conditions of $H_n$ will not
be known explicitly in general. To evade this obstacle they show that their
result still holds if the Weyl solution of a nearby operator is chosen instead.

The main purpose of this note is to propose an alternate way of approximating $H$
which only involves one Weyl solution at one endpoint. More precisely, we
show that if $\tau$ is l.p.\ at one end point, the Weyl
solution of the other endpoint can also be chosen instead:

\begin{theorem}\label{thm:main}
Suppose $\tau$ is l.p.\ at $b$.
Define $H_n$ as above with $u=\psi_a(\lam_a)$, $\lam_a\in[\lam_0,\lam_1]$, and
$v=\psi_a(\lam_0)$ or $v=\psi_a(\lam_1)$
(in particular, we assume that the corresponding Weyl solutions exist).

If $\dim\Ran P_{(\lam_0,\lam_1)}(H)<\infty$, then the eigenvalues of $H$
in $(\lam_0,\lam_1)$ are exactly the limits of eigenvalues of $H_n$ which lie in
$(\lam_0,\lam_1)$. The corresponding (one-dimensional) eigenprojections converge
in norm.

If $\dim\Ran P_{(\lam_0,\lam_1)}(H)=\infty$, then the eigenvalues of $H_n$
accumulate in $(\lam_0,\lam_1)$ as $n\to\infty$.
\end{theorem}

\noindent
In particular, if one endpoint is regular, a solution satisfying the boundary condition
at this endpoint can be taken in this case. Clearly the result cannot hold if $\tau$
is l.c.\ at $b$, since our assumptions contain no information on the boundary condition
of $H$ at $b$ in this case.

\begin{remark}
(i). As shown in \cite{sw}, if $\tau$ is l.c.\ at $a$, then $u$ can be choosen to be any
function in $\db(\tau)$ generating the boundary condition of $H$ at $a$.

(ii). The same result (with the same proof) holds for Jacobi operators (see \cite{tjac}).
\end{remark}

\section{Approximation by regular operators}

We begin by recalling that $H_n$ converges to $H$ strongly (\cite{wdln}). Strictly speaking
this statement makes no sense  since $H_n$ and $H$ live in different Hilbert spaces.
This can be easily fixed by using $\alpha \id \oplus H_n \oplus\alpha \id$ on
$L^2((a,b);r\,dx)=L^2((a,b_n);r\,dx)\oplus L^2((a_n,b_n);r\,dx)\oplus L^2((a_n, b);r\,dx)$ where
$\alpha$ is a fixed real constant outside $[\lam_0,\lam_1]$. Alternatively, one can
also use generalized strong convergence as introduced in \cite{sw}.

\begin{lemma} \label{strconh}
Suppose that either $H$ is limit point at $a$ or
that $u = \psi_-(\lam_0)$ for some $\lam_0$ and similarly, that either
$H$ is limit point at $b$ or $v=\psi_+(\lam_1)$ for some $\lam_1$.
Then $H_m$ converges to $H$ in strong resolvent sense as $m\to\infty$.
\end{lemma}

\noindent
In addition, we need the following abstract result. In this respect
we remark that for a self-adjoint projector $P$ we have
\be
\dim\Ran(P) = \tr(P) = \|P\|_1,
\ee
where $\|.\|_1$ denotes the trace class norm. If $P$ is not finite rank, then
it is of course not trace class and all three numbers are equal $\infty$.

\begin{lemma}
Let $A_n$, $A$ be self-adjoint operators such that $A_n\to A$ in strong resolvent sense.
Then
\be
\tr(P_{(\lam_0,\lam_1)}(A)) \le \liminf \tr(P_{(\lam_0,\lam_1)}(A_n)).
\ee
If in addition, $\tr(P_{(\lam_0,\lam_1)}(A_n)) \le \tr(P_{(\lam_0,\lam_1)}(A))$, then
\be
\lim_{n\to\infty} \tr(P_{(\lam_0,\lam_1)}(A_n)) = \tr(P_{(\lam_0,\lam_1)}(A))
\ee
and if $\tr(P_{(\lam_0,\lam_1)}(A))<\infty$ we even have
\be
\lim_{n\to\infty} \| P_{(\lam_0,\lam_1)}(A_n) - P_{(\lam_0,\lam_1)}(A)\|_1 =0.
\ee
\end{lemma}

\begin{proof}
The first part is just Lemma~5.2 from \cite{gst}. This also implies the second if
$\tr(P_{(\lam_0,\lam_1)}(A))=\infty$. Otherwise, if $\tr(P_{(\lam_0,\lam_1)}(A))<\infty$,
we have
$$
\limsup \tr(P_{(\lam_0,\lam_1)}(A_n)) \le \tr(P_{(\lam_0,\lam_1)}(A))
$$
and the first claim follows. The second is then a consequence of
Gr\"umm's theorem (\cite[Thm.~2.19]{str}).
\end{proof}

\noindent
Now it remains to show that this result is applicable in our situation.

\begin{lemma}
Let $H_n$ be defined as in (\ref{deftiHm}) with 
\begin{enumerate}
\item
$u = \psi_a(\lam_a)$, $v=\psi_b(\lam_b)$ with $\lam_a,\lam_b \in[\lam_0,\lam_1]$ or
\item
$u = \psi_a(\lam_a)$ with $\lam_a \in [\lam_0,\lam_1]$ and
$v = \psi_a(\lam_b)$ with $\lam_b \in \{\lam_0,\lam_1\}$.
\end{enumerate}
Then, 
\be
\tr(P_{(\lam_0, \lam_1)}(H_n)) \le  \tr(P_{(\lam_0, \lam_1)}(H)).
\ee
\end{lemma}

\begin{proof}
Abbreviate $P=P_{(\lam_0, \lam_1)}(H)$, $P_n=P_{(\lam_0, \lam_1)}(H_n)$.

(i). Since this part is identical to the proof in \cite{sw}, we just give an outline.
Let $\ti{\psi}_1, \dots, \ti{\psi}_k \in \Ran P_n$ be the normalized eigenfunctions of $H_n$,
construct
$$
\psi_j (x) = \left\{\begin{array}{cl} \gam_{a,j} u(x), & x < a_n,\\
\ti{\psi}_j(x), & a_n \leq x \leq b_n,\\
\gam_{b,j} v(x), & x > b_n, \end{array} \right.
$$
where $\gam_{a,j}$, $\gam_{b,j}$ are chosen such that $\psi_j\in\db(\tau)$.
A computation now shows that
$$
\|(H - \frac{\lam_1 + \lam_0}{2}) \psi\| < \frac{\lam_1 - \lam_0}{2} \|\psi\|
$$
for any $\psi$ in the linear span of the $\psi_j$'s, which yields the first result.

(ii) By considering two steps from $(a_n,b_n)$ to $(a,b_n)$ and from 
$(a,b_n)$ to $(a,b)$, we see that the first step is covered by (i) and hence
it is no restriction to assume $a_n=a$. Now proceed as in the previous case
but use
$$
\psi_j (x) = \left\{\begin{array}{cl}
\ti{\psi}_j(x) - \gam_j v(x) , & x \leq b_n,\\
0, & x > b_n, \end{array} \right.
$$
where $\gam_j$ are chosen such that $\psi_j\in\db(\tau)$. Now
let $\psi =\sum_j c_j \psi_j$ be in the linear span of the $\psi_j$'s.
Then, since $v$ is also an eigenvector of $H_n$ and hence orthogonal to
the $\psi_j$'s, we have
\begin{align*}
\|(H - \frac{\lam_1 + \lam_0}{2}) \psi\|^2 &=
\| \sum_j c_j \Big( \frac{2\lam_j - \lam_1-\lam_0}{2}\ti{\psi}_j(x) -
\frac{2\lam_b - \lam_1-\lam_0}{2} \gam_j v(x) \Big)\|^2\\
&= \sum_j |c_j|^2 (\frac{2\lam_j - \lam_1-\lam_0}{2})^2 +
|\gam|^2 (\frac{\lam_1-\lam_0}{2})^2 \|v\|^2_{(a,b_n)}\\
&< (\frac{\lam_1-\lam_0}{2})^2 \|\psi\|^2,
\end{align*}
where $\gam=\sum_j c_j \gam_j$. Hence the second result follows.
\end{proof}

\noindent
Theorem~\ref{thm:sw} and Theorem~\ref{thm:main} now follow by combining the
last two lemmas.

\section*{Acknowledgments}

I thank Helge Kr\"uger for discussions.

\end{document}